\documentclass[preprint,12pt]{elsarticle}

\usepackage{graphicx}
\usepackage{amssymb}

\usepackage{showlabels}
\usepackage{epsf}
\usepackage{amsbsy,amsmath}
\usepackage{mathtools}
\usepackage{mathrsfs}
\usepackage{amsfonts}
\usepackage{amssymb}
\usepackage{enumitem}
\usepackage{eucal}
\usepackage{cases}
\usepackage{graphics,mathrsfs}
\usepackage{amsthm}
\usepackage{secdot}
\usepackage{esint}
\usepackage{hyperref}
\usepackage{varwidth}
\usepackage{tasks}
\usepackage{geometry}
\newtheorem{theorem}{Theorem}[section]
\newtheorem{lemma}[theorem]{Lemma}
\newtheorem{proposition}[theorem]{Proposition}

\theoremstyle{definition}

\theoremstyle{remark}
\newtheorem{remark}[theorem]{Remark}
\numberwithin{equation}{section}

\numberwithin{equation}{section}

\usepackage{xcolor}

\begin{document}



 \begin{frontmatter}
 	
 	
 	
 	\title{A singular elliptic problem involving fractional $p$-Laplacian and a discontinuous critical nonlinearity}
 	 \author{K. Saoudi$^{\dagger}$\corref{cor1}\fnref{label3}}
 	\ead{kasaoudi@gmail.com}
  \author{Akasmika Panda\fnref{label2}}
  \ead{akasmika44@gmail.com}
  	\author{Debajyoti Choudhuri\fnref{label2}}
  \ead{dc.iit12@gmail.com}

 	%
 	\cortext[cor1]{Corresponding author - K. Saoudi}
 	\fntext[label2]{Department of Mathematics, National Institute of Technology Rourkela}
\fntext[label3]{Basic and Applied Scientific Research Center, Imam Abdulrahman Bin Faisal University,  Saudi Arabia}


\begin{abstract}
	In this article, we prove the existence of solutions to a  nonlinear nonlocal elliptic problem with a singualrity and a discontinuous critical nonlinearity which is given as follows.
	\begin{align}
		\begin{split}\label{main_prob}
			(-\Delta)_p^su&=\mu g(x,u)+\frac{\lambda}{u^\gamma}+H(u-\alpha)u^{p_s^*-1},~\text{in}~\Omega\\
		    u&>0,~\text{in}~\Omega,\\
			u&=0,~\text{in}~\mathbb{R}^N\setminus\Omega,
		\end{split}
	\end{align}
	where $\Omega\subset\mathbb{R}^N$ is a bounded domain with Lipschitz boundary, $s\in (0,1)$, $2<p<\frac{N}{s}$, $\gamma\in (0,1)$, $\lambda,\mu>0$, $\alpha\geq 0$ is real, $H$ is the Heaviside function, i.e. $H(a)=0$ if $a\leq 0$, $H(a)=1$ if $a>0$ and $p_s^*=\frac{Np}{N-sp}$ is the fractional critical Sobolev exponent.\\
	Under suitable assumptions on the function $g$, we prove the existence of solution to the problem. Furthermore, we show that as $\alpha\rightarrow0^+$, the sequence of solutions of $\eqref{main_prob}$ for each such $\alpha$ converges to a solution of the problem for which $\alpha=0$.
	\begin{flushleft}
		{\bf Keywords}:~  Fractional  $p$-Laplacian, Heaviside function, Mountain pass theorem, Critical exponent, Singularity.\\
		{\bf AMS Classification}:~35R11, 35J75, 35J60, 46E35.
	\end{flushleft}
\end{abstract}
\end{frontmatter}

\section{Introduction}\label{intro}
\noindent We will study the existence of solution to the following nonlinear, nonlocal problem involving a singularity and a discontinuous critical nonlinearity. 
	\begin{align}\tag{$P_\alpha$}
	\begin{split}\label{main prob intro}
	(-\Delta)_p^su&=\mu g(x,u)+\frac{\lambda}{u^\gamma}+H(u-\alpha)u^{p_s^*-1},~\text{in}~\Omega\\
	u&>0,~\text{in}~\Omega,\\
	u&=0,~\text{in}~\mathbb{R}^N\setminus\Omega,
	\end{split}
	\end{align}
	We impose the hypotheses on $g$ which are as follows.
	\begin{itemize}
		\item[($g_1$)] the function $g:\Omega\times\mathbb{R}\rightarrow\mathbb{R}$ is a $N$-measurable function and $g(x,u)=0$ if $u\leq 0$.
		\item[($g_2$)] there exist $K>0$ and $r\in (p,p_s^*)$ such that $|g(x,u)|\leq K(1+|u|^{r-1})$ for every $u\geq 0$.
		\item[($g_3$)] there exists $b>p$ and $v>0$ such that for all $u\geq v$ $$0<bG(x,u)=b\int_{0}^{u}g(x,\tau)d\tau\leq u\underline{g}(x,u).~\text{refer Section}~\ref{priliminaries}~\text{for the definition of $\underline{g}$}$$
		\item[($g_4$)] there exists $\beta>0$ (that will be fixed later) such that $H(u-\beta)\leq g(x,u)$ uniformly in $\Omega\times(0,\infty)$.
		\item[($g_5$)] Let $\lambda_1$ be the first eigen value of $(-\Delta)_p^s$ defined in $\eqref{first}$. Then 
		$\lim\limits_{u\rightarrow0}\frac{g(x,u)}{u^{p-1}}\leq \lambda_1$ uniformly in $\Omega$.
	\end{itemize}
	A prototype of $g$ satisfying the asuumptions $(g_1)-(g_5)$ is $H(t-\beta)t^{r-1}/\beta^r$.\\
	The problems of type $\eqref{main prob intro}$ having discontinuous nonlinearities have many applications in free boundary problems of mathematical physics. For instance, obstacle problem, Elenbaas equations, the seepage surface problem etc. Refer \cite{Ambrosetti,Chang,Chang 2,Chang 3} for further details.\\
	Elliptic problems involving critical and discontinuous nonlinearities can be treated by different techniques. Amongst these methods, variational methods for nondifferentiable functionals, dual variational principle, Palais principle of symmetric criticality for locally Lipschitz functional, lower-upper solution method, theory of multivalued mappings and global branching are a few well known techniques. Badiale \& Tarantello in \cite{Badiale 2} studied the following class of problem using variational methods with lower-upper solution methods.
		\begin{align}\label{prob 1}
		\begin{split}
		(-\Delta)u&=\delta H(u-\alpha)+u^{2_s^*-1},~\text{in}~\Omega\\
		u&=0,~\text{on}~\partial \Omega.
		\end{split}
		\end{align}
Here $2^*=2N/(N-2)$, $\delta,\alpha>0$ and $H$ is the Heaviside function. Later, the authors in \cite{Alves 2} generalized the work of \cite{Badiale 2} in $\mathbb{R}^N$. Badiale in \cite{Badiale 1} proved the existence result for the critical elliptic problem given by
	\begin{align}\label{prob 2}
	\begin{split}
	(-\Delta)u&=g(u)+u^{2_s^*-1},~\text{in}~\Omega\\
	u&=0,~\text{on}~\partial \Omega,
	\end{split}
	\end{align}
	where $g$ can have discontinuities. The authors in \cite{Santos} and \cite{Figueiredo}
 extended the result of \cite{Badiale 1} for a Kirchhoff type problem involving critical Caffarelli-Kohn-Nirenberg growth and for a Schr\"{o}dinger-Kirchhoff equation, respectively. Recently, Santos \& Tavares in \cite{Santos 2} considered the problem
 \begin{align}
 \begin{split}\label{prob 3}
 \mathcal{L}_Ku&= g(x,u)+H(u-\alpha)u^{2_s^*-2}u,~\text{in}~\Omega\\
 u&=0,~\text{in}~\mathbb{R}^N\setminus\Omega,\\
 u&\geq 0,~\text{in}~\Omega,
 \end{split}
 \end{align}
 where $2_s^*=2N/(N-2s)$, $\alpha> 0$, $g$ is a discontinuous function and $\mathcal{L}_K$ is a nonlocal operator 
 $$\mathcal{L}_K u(x)=\iint_{\mathbb{R}^{2N}}(u(x+y)+u(x-y)-2u(x))K(y)dy.$$
 They used the nonsmooth version of Mountain pass theorem to investigate the existence and the behavior of solution for problem $\eqref{prob 3}$. We also cite \cite{Alves 1,Ambrosetti,Arcoya,Bananno,Chang,Chang 3,Clarke,Radulescu,Xiang} and the references therein for readers to have a glimpse of the problems of the type as in $\eqref{prob 1}-\eqref{prob 3}$.\\
 Inspired by the above works, specifically \cite{Badiale 1,Badiale 2,Chang,Chang 2,Santos 2} we analyze our problem $\eqref{main prob intro}$. The problem $\eqref{main prob intro}$ with singularity, critical and discontinuous nonlinearities is a new and first work in the literature, at least to our knowledge. But we find enormous works dealing with the following class of problems involving singularity and critical nonlinearity given by
 \begin{align}
 \begin{split}\label{prob 4}
 (-\Delta)_p^su&=\frac{\lambda g_1(x)}{u^\gamma}+\delta u^{q-1},~\text{in}~\Omega\\
 u&=0,~\text{in}~\mathbb{R}^N\setminus\Omega,\\
 u&>0,~\text{in}~\Omega,
 \end{split}
 \end{align}
 where $\lambda,\delta>0$, $q\in (1,p_s^*]$, $g_1>0$ is bounded. Several techniques like variational method, concentration compactness method, Nehari manifold method etc. have been applied to study the problems of type $\eqref{prob 4}$ for both local and nonlocal cases. Refer \cite{Dhanya,Ghanmi,Ghosh,Giacomoni,Giacomoni 1,Giacomoni 2,Haitao,Hirano,Mukherjee,Mukherjee 1} and the bibliography therein.\\	
The main result of this article is the following.
	\begin{theorem}\label{main result heaviside}
	Let $(g_1)-(g_5)$ hold. Then 
	\begin{enumerate}
		\item there exist $\bar{\alpha},\bar{\lambda}, \bar{\mu}>0$ such that for every $a\in (0,\bar{\alpha})$, every $\lambda\in (0,\bar{\lambda})$ and every $\mu\in (0,\bar{\mu})$ the problem $\eqref{main prob intro}$ admits at least one nontrivial weak solution $u_\alpha$. Furthermore, the lebesgue measure of the set $\{x\in\Omega:u_\alpha>\alpha\}$ is positive
		\item for any sequence $\alpha_n\in (0,\bar{\alpha})$ with $\alpha_n\rightarrow 0^+$, we have, up to a subsequential level, $u_{\alpha_n}\rightarrow u_0$ in $W_0^{s,p}(\Omega)$, where $u_0$ is a weak solution of the problem $(P_0)$, i.e. 
		\begin{align}\tag{$P_0$}
		\begin{split}\label{main prob 0}
		(-\Delta)_p^su&=\mu g(x,u)+\frac{\lambda}{u^\gamma}+u^{p_s^*-1},~\text{in}~\Omega\\
		u&=0,~\text{in}~\mathbb{R}^N\setminus\Omega,\\
		u&>0,~\text{in}~\Omega.
		\end{split}
		\end{align}
	\end{enumerate}	
	\end{theorem}
	\noindent The proof of the main result, i.e Theorem $\ref{main result heaviside}$, has been splitted into two sections. Section $\ref{3}$ is devoted to the first part of Theorem $\ref{main result heaviside}$, i.e. the existence of a weak solution $u_\alpha$ to $\eqref{main prob intro}$. In Section $\ref{4}$, we examine the nature of the sequence $(u_\alpha)$ as $\alpha\rightarrow0^+$ and prove the second part of Theorem $\ref{main result heaviside}$.
	\section{Nonsmooth critical point theory}\label{priliminaries}
\noindent Let us fix $0<s<1$, $2<p<\frac{N}{s}$, $\Omega$ be an open and bounded domain of $\mathbb{R}^N$. We denote $Q=(\mathbb{R}^N\times\mathbb{R}^N)\setminus(\Omega^c\times\Omega^c)$ where $\Omega^c=\mathbb{R}^N\setminus\Omega$. We define the fractional Sobolev space  by
	\begin{align} 
	W^{s,p}(\Omega)=\{u:\mathbb{R}^N\rightarrow\mathbb{R}~\text{is measurable}:u|_\Omega\in L^{p}(\Omega),\int_{Q}\frac{|u(x)-u(y)|^{p}}{|x-y|^{N+sp}}dydx<\infty\}\nonumber
	\end{align}
	 equipped with the norm 
	$$\|u\|_{W^{s,p}(\Omega)}=\|u\|_{L^{p}(\Omega)}+\left(\int_{Q}\frac{|u(x)-u(y)|^{p}}{|x-y|^{N+sp}}dydx\right)^{1/p}.$$
	We further define the space
	$$W_0^{s,p}(\Omega)=\{u\in 	W^{s,p}(\Omega):u=0 ~\text{a.e. in}~\mathbb{R}^N\setminus\Omega\}$$ 
and $(W_0^{s,p}(\Omega),\|\cdot\|_{W_0^{s,p}(\Omega)})$ is a reflexive Banach space where the fractionl Sobolev norm is given by
	$$\|u\|^p_{W_0^{s,p}(\Omega)}=\int_{Q}\frac{|u(x)-u(y)|^{p}}{|x-y|^{N+sp}}dydx.$$
	 Given below are a few well known embedding results for the space $W_0^{s,p}(\Omega)$.
	 \begin{theorem}[\cite{Servadei}]\label{embedding}
	 	The following results hold for the fractional Sobolev space $W_0^{s,p}(\Omega)$.
	 	\begin{enumerate}
	 		\item If $\Omega$ has a continuous boundary, then the embedding $W_0^{s,p}(\Omega)\hookrightarrow L^q(\Omega)$ is compact for every $q\in[1,p_s^*)$.
	 		\item The embedding $W_0^{s,p}(\Omega)\hookrightarrow L^{p_s^*}(\Omega)$ is continuous.
	 	\end{enumerate}
	 \end{theorem}
	\noindent We now define the best constant $S_{s,p}>0$ given by
	\begin{equation}\label{best}
	S_{s,p}=\underset{u\in W_0^{s,p}(\Omega)\setminus\{0\}}{\inf}\frac{\int_{Q}\frac{|u(x)-u(y)|^{p}}{|x-y|^{N+sp}}dydx}{(\int_\Omega |u|^{p_s^*}dx)^{\frac{p}{p_s^*}}}
	\end{equation}
and $S_{s,p}$ is well-defined due to Theorem $\ref{embedding}$.
\begin{theorem}[\cite{Brasco}, Thoerem 4.1 of \cite{Franzina}]
	Let $s\in (0,1)$ and $p>1$. Then the eigenvalue problem
	\begin{align}
	\begin{split}
	(-\Delta)_p^su&=\lambda |u|^{p-2}u,~\text{in}~\Omega\\
	u&=0,~\text{in}~\mathbb{R}^N\setminus\Omega\nonumber
	\end{split}
	\end{align}
 possesses a smallest eigenvalue $\lambda_1>0$ given by
 \begin{equation}\label{first}
 \lambda_1=\underset{u\in W_0^{s,p}(\Omega):\|u\|_{L^{p}(\Omega)}=1}{\min}\|u\|_{W_0^{s,p}(\Omega)}^p. 
 \end{equation}
\end{theorem}	
\noindent	Let $V$ be a Banach space. Then a functional $J:V\rightarrow\mathbb{R}$ is said to be locally Lipschitz continuous if for any $u\in V$ there exists an open neighborhood $N:=N_u\subset V$ and some constant $C:=C_N>0$ such that
$$|J(u_1)-J(u_2)|\leq C\|u_1=u_2\|_V, ~u_1,u_2\in V.$$
We define the directional derivative of $J$ at $u$ in the direction of $z\in V$ by
$$\tilde{J}(u;z)=\lim\limits_{h\rightarrow0}\lim\limits_{\xi\downarrow0}\frac{J(u+h+\xi z)-J(u+h)}{\xi}.$$
Thus, $\tilde{J}(u;\cdot)$ is convex, continuous and its subdifferential at $w\in V$ is the set
$$\partial\tilde{J}(u;w)=\{\nu\in V^*:\tilde{J}(u;z)\geq\tilde{J}(u;w)+\langle \nu, z-w\rangle, ~z\in V\}.$$
Here $\langle\cdot,\cdot\rangle$ is the duality pair between $V$ and $V^*$. The generalized gradient of $J$ at $u$ is defined by
$$\partial J(u)=\{\nu\in V^*:\langle \nu, z\rangle\leq \tilde{J}(u;z), ~z\in V\}$$
and is convex and weak$^*$- compact. Clearly, $\partial J(u)$ is nonempty and is the subdifferential of $\partial\tilde{J}(u;0)$ as $\tilde{J}(u;0)=0$.\\
We say $\bar{u}$ to be a critical point of $J$ if $0\in \partial J(\bar{u})$ and $c\in 
\mathbb{R}$ to be a critical value of $J$ if $J(\bar{u})=c$ for some critical point $\bar{u}\in V$. Moreover, if $J$ is a $C^1$ functonal then $\partial J(u)=\{J^\prime(u)\}$. We now denote
$$\Lambda_J(u)=\min\{\|\nu\|_{V^*}:\nu\in \partial J(u)\}.$$
The following result is the Mountain Pass Theorem for locally Lipschitz non-differentiable functional.
\begin{theorem}[\cite{Grossinho,Radulescu}]\label{mountain pass non smooth}
Let $V$ be a Banach space and $J$ be a locally Lipschitz functional with $J(0)=0$. Assume that there exists $\rho_1,\rho_2>0$ and $\sigma\in V$ such that 
\begin{enumerate}
	\item $J(u)\geq \rho_2$, for every $u\in V$; $\|u\|_V=\rho_1$,
	\item $J(\sigma)<0$ and $\|\sigma\|_{V}>\rho_1$.
\end{enumerate}
Let $$c=\underset{\zeta\in \Gamma}{\inf}\underset{t\in [0,1]}{\max}~J(\zeta(t)),~\Gamma=\{\zeta\in C([0,1];V):\zeta(0)=0,\zeta(1)=\sigma\}.$$
Then $c\geq \rho_2$ and there exists a Palais-Smale [(PS)$_c$] sequence $(u_n)\subset V$ satisfying $$J(u_n)\rightarrow c~\text{and}~\Lambda_J(u_n)\rightarrow 0.$$
Moreover, if $J$ satisfies the nonsmooth Palais-Smale [(PS)$_c$] condition, i.e. every $(PS)_c$ sequence has a convergent subsequence, then $c$ is a critical value of $J$. 
\end{theorem}
\begin{proposition}[\cite{Grossinho}]\label{convergence}
	Let $(u_n)\subset V $ and $(\theta_n)\subset V^*$ with $\theta_n\in \partial J(u_n)$. If $u_n\rightarrow u$ in $V$ and $\theta_n\overset{*}{\rightharpoonup}\theta$ in $V^*$, then $\theta\in \partial J(u)$.
\end{proposition}
\noindent A function $g:\Omega\times\mathbb{R}\rightarrow\mathbb{R}$ is said to be a $N$-measurable function if for every $u\in W_0^{s,p}(\Omega)$ the function $g(\cdot,u(\cdot)):V\rightarrow \mathbb{R}$ is measurable (refer \cite{Chang}). Let $g(x,\cdot)\in L_{loc}^1(\Omega)$, then we denote
$$\underline{g}(x,u)=\underset{\delta\downarrow0}{\lim}~\underset{|v-u|<\delta}{essinf}~g(x,v),~\overline{g}(x,u)=\underset{\delta\downarrow0}{\lim}~\underset{|v-u|<\delta}{esssup}~g(x,v).$$
We now provide the notion of weak solution to $\eqref{main prob intro}$ (influenced by \cite{Santos 2}). We say a function $u\in W_0^{s,p}(\Omega)$ to be a weak solution of $\eqref{main prob intro}$ if $u>0$ a.e. in $\Omega$, $u^{-\gamma}\in L^1_{loc}(\Omega)$ and there exist $\eta_\alpha\in L^{\frac{r}{r-1}}(\Omega)$ and $\theta_\alpha\in L^{\frac{p_s^*}{p_s^*-1}}(\Omega)$ such that for every $\varphi\in W_0^{s,p}(\Omega)$
\begin{equation}
\int_{\mathbb{R}^{2N}}\frac{|u(x)-u(y)|^{p-2}(u(x)-u(y)}{|x-y|^{N+sp}}(\varphi(x)-\varphi(y))dxdy-\mu\int_\Omega \eta_\alpha \varphi dx-\int_\Omega\theta_\alpha\varphi dx-\lambda\int_\Omega\frac{\varphi}{u^{\gamma}}dx=0.
\end{equation}
Here $\eta_\alpha(x)\in [\underline{g}(x,u(x)),\overline{g}(x,u(x))]$ and $\theta_\alpha(x)\in [\underline{f}_\alpha(u(x)),\overline{f}_\alpha(u(x))]$ a.e. in $\Omega$ with $f_\alpha(t)=H(t-\alpha)t^{p_s^*-1}$.
\section{Existence result for $\eqref{main prob intro}$ - Proof of Theorem $\ref{main result heaviside}$ (1)}\label{3}
\noindent The energy functional corresponding to the problem $\eqref{main prob intro}$ is $J_\alpha^0:W_0^{s,p}(\Omega)\rightarrow\mathbb{R}$ defined by
$$J_\alpha^0(u)=\frac{1}{p}\int_{\mathbb{R}^{2N}}\frac{|u(x)-u(y)|^{p}}{|x-y|^{N+sp}}dxdy-\frac{\lambda}{1-\gamma}\int_\Omega (u^+)^{1-\gamma}-\mu\int_\Omega G(x,u)-\int_\Omega F_\alpha(u), $$
where $G(x,u)=\int_{0}^{u}g(x,\tau)d\tau$, $ F_\alpha(u)=\int_{0}^{u}f_\alpha(\tau)d\tau$ with $f_\alpha(u)=H(u-\alpha)u^{p_s^*-1}$. Let us denote
\begin{equation}\label{I, I alpha}
I_\lambda^0(u)= \frac{1}{p}\int_{\mathbb{R}^{2N}}\frac{|u(x)-u(y)|^{p}}{|x-y|^{N+sp}}dxdy-\frac{\lambda}{1-\gamma}\int_\Omega (u^+)^{1-\gamma},~I_\mu(u)=\int_\Omega G(x,u),~ I_\alpha(u)=\int_\Omega F_\alpha(u).
\end{equation} 
We provide some properties of the functionals $I_\mu$ and $I_\alpha$ in the following lemma and these results can be proved by following the arguments of Lemma 3.1 of \cite{Santos}.
\begin{lemma}\label{properties I}
	The functionals $I_\mu:L^r(\Omega)\rightarrow\mathbb{R}$, $I_\alpha:L^{p_s^*}(\Omega)\rightarrow\mathbb{R}$ are locally Lipschitz functionals and they satisfy the followings.
	\begin{enumerate}
		\item $\partial I_\mu(u)\subset [\underline{g}(x,u(x)),\overline{g}(x,u(x))]$ a.e. in $\Omega$.
		\item $\partial I_\alpha(u)\subset [\underline{f}_\alpha(u(x)),\overline{f}_\alpha(u(x))]$ a.e. in $\Omega$.
		\item \begin{eqnarray}
		\begin{split}
		[\underline{f}_\alpha(u),\overline{f}_\alpha(u)]=&\begin{cases}
		\{0\}, & ~\text{if}~u<a\\ [0,u^{p_s^*-1}], &~\text{if}~u=a\\\{u^{p_s^*-1}\}, &~\text{if}~u>a.
		\end{cases}
		\end{split}
		\end{eqnarray}
	\end{enumerate}
\end{lemma}
\begin{remark}
The inclusion (1) of Lemma $\ref{properties I}$ imply that if $\eta\in \partial I_\mu(u)$ then $\eta\in L^{\frac{r}{r-1}}(\Omega)$ and $\eta(x)\subset [\underline{g}(x,u(x)),\overline{g}(x,u(x))]$ a.e. in $\Omega$. The same argument follows for $(2)$.
\end{remark}
\noindent The functional $I_\lambda^0$ is not of class $C^1$ due to the presesnce of a singularity. To tackle this issue, we follow the truncation technique. Let us consider the problem
	\begin{align}
	\begin{split}\label{subsolution}
	(-\Delta)_p^su&=\frac{\lambda}{u^\gamma},~\text{in}~\Omega\\
	u&=0,~\text{in}~\mathbb{R}^N\setminus\Omega,\\
	u&>0,~\text{in}~\Omega.
	\end{split}
	\end{align}
	\noindent According to Canino et al. in \cite{Canino}, we have the existence result for $\eqref{subsolution}$ as given below.
	\begin{lemma}\label{main result help}
		Let $\gamma\in(0,1)$ and $\lambda>0$. Then $\eqref{subsolution}$ admits a unique nontrivial weak
		solution $\underline{u}_\lambda$ in $W_0^{s,p}(\Omega)$ such that for every $\omega\subset\subset\Omega$  we have $\underset{\omega}{ess~inf}~\underline{u}_\lambda>0$.
	\end{lemma}
	\noindent We now define
	\begin{eqnarray}\label{defn}
	\begin{split}
	\psi(x,t)=&\begin{cases}
	t^{-\gamma}, & ~\text{if}~t>\underline{u}_\lambda\\ \underline{u}_\lambda^{-\gamma}, &~\text{if}~t\leq \underline{u}_\lambda,
	\end{cases}
	\end{split}
	\end{eqnarray}
	where $\underline{u}_\lambda$ is the unique solution to $\eqref{subsolution}$. Further, define a function
 $I_\lambda:W_0^{s,p}(\Omega)\rightarrow\mathbb{R}$ by
 	\begin{equation}\label{I lambda}
 	I_\lambda(u)=\frac{1}{p}\int_{\mathbb{R}^{2N}}\frac{|u(x)-u(y)|^{p}}{|x-y|^{N+sp}}dxdy-\lambda\int_\Omega \Psi(x,u)dx
 	\end{equation}
where $\Psi(x,t)=\int_{0}^{t}\psi(x,\tau)d\tau$.\\
 We now consider the following cutoff problem.
	\begin{align}
	\begin{split}\label{main prob help}
	(-\Delta)_p^sw&=\mu g(x,w)+\lambda \psi(x,w)+H(w-\alpha)w^{p_s^*-1},~\text{in}~\Omega\\
	w&=0,~\text{in}~\mathbb{R}^N\setminus\Omega,\\
	w&>0,~\text{in}~\Omega.
	\end{split}
	\end{align}
	 A function $w\in W_0^{s,p}(\Omega)$ is said to be a weak solution of $\eqref{main prob help}$ if $w>0$ a.e. in $\Omega$ and there exist $\eta_\alpha\in L^{\frac{r}{r-1}}(\Omega)$ and $\theta_\alpha\in L^{\frac{p_s^*}{p_s^*-1}}(\Omega)$ such that for every $\varphi\in W_0^{s,p}(\Omega)$
	 \begin{equation}
	 \int_{\mathbb{R}^{2N}}\frac{|w(x)-w(y)|^{p-2}(w(x)-w(y)}{|x-y|^{N+sp}}(\varphi(x)-\varphi(y))dxdy-\mu\int_\Omega \eta_\alpha \varphi dx-\int_\Omega\theta_\alpha\varphi dx-\lambda\int_\Omega \psi(x,w) \varphi dx=0,
	 \end{equation}
	 where $\eta_\alpha(x)\in [\underline{g}(x,w(x)),\overline{g}(x,w(x))]$ and $\theta_\alpha(x)\in [\underline{f}_\alpha(w(x)),\overline{f}_\alpha(w(x))]$ a.e. in $\Omega$.\\
	 The associated functional of $\eqref{main prob help}$ is $J_\alpha:W_0^{s,p}(\Omega)\rightarrow\mathbb{R}$ defined by
	 \begin{equation}\label{J alpha}
	 J_\alpha(w)=I_\lambda(w)-I_\mu(w)-I_\alpha(w).
	 \end{equation}
The functionals $I_\lambda, I_\mu$ and $I_\alpha$ are given in $\eqref{I, I alpha}$ and $\eqref{I lambda}$. By Lemma $\ref{properties I}$, the functional $J_\alpha$ is a locally Lipschitz functional on $W_0^{s,p}(\Omega)$ and thus using Proposition 1.3.12 and 1.3.13 of \cite{Gasinski} we obtain 
\begin{equation}\label{differential}
\partial J_\alpha(w)\subset \{I_\lambda^\prime(w)\}-\partial I_\mu(w)-\partial I_\alpha(w), ~\text{for all}~w\in W_0^{s,p}(\Omega).
\end{equation}
It is easy to check that if $w$ is a weak solution to $\eqref{main prob help}$ with $w\geq \underline{u}_\lambda$ a.e. in $\Omega$, then $w$ is also a weak solution to $\eqref{main prob intro}$. Hence, with the help of a comparison principle and non-smooth variational approach we prove our main result.
\begin{lemma}\label{lemma 1}
	Any sequence $(w_n)\subset W_0^{s,p}(\Omega)$ satisfying $J_\alpha(w_n)\rightarrow c$ and $\Lambda_{J_\alpha}(w_n)\rightarrow 0$ is bounded in $W_0^{s,p}(\Omega)$.
	\begin{proof}
		Suppose $(v_n)\subset (W_0^{s,p}(\Omega))^\prime$ is a sequence with $\|v_n\|_{(W_0^{s,p}(\Omega))^\prime}=\Lambda_{J_\alpha}(w_n)$, where $(W_0^{s,p}(\Omega))^\prime$ is the dual of $W_0^{s,p}(\Omega)$. This implies $v_n\in \partial J_\alpha(w_n)$. From $\eqref{differential}$, there exist $\eta_n\in \partial I_\mu(w_n)$ and $\theta_n\in \partial I_\alpha(w_n)$ that satisfy
		\begin{align}\label{representation}
		\langle v_n,\varphi\rangle&=\int_{\mathbb{R}^{2N}}\frac{|w_n(x)-w_n(y)|^{p-2}(w_n(x)-w_n(y)}{|x-y|^{N+sp}}(\varphi(x)-\varphi(y))dxdy-\lambda\int_\Omega \psi(x,w_n)\varphi dx\nonumber\\&~~~~-\mu\langle\eta_n,\varphi\rangle-\langle\theta_n,\varphi\rangle,~\forall \varphi\in W_0^{s,p}(\Omega).
		\end{align}
		By $\eqref{J alpha}$, $\eqref{representation}$ and Lemma $\ref{properties I}$, we get 
		\begin{align}\label{star}
		J_\alpha(w_n)-\frac{1}{b}\langle v_n,w_n\rangle=&\left(\frac{1}{p}-\frac{1}{b}\right)\|w_n\|^p_{W_0^{s,p}(\Omega)}+\frac{\lambda}{b}\int_\Omega \psi(x,w_n)w_n-\lambda\int_\Omega \Psi(x,w_n)\nonumber\\&~~~~+\mu\int_\Omega \left(\frac{1}{b}\eta_n w_n-G(x,w_n)\right)+\int_\Omega \left(\frac{1}{b}\theta_n w_n-F_\alpha(w_n)\right).
		\end{align}
		According to $(g_3)$,
	\begin{equation}\label{star 1}
	\int_\Omega \left(\frac{1}{b}\eta_n w_n-G(x,w_n)\right)\geq\int_{\{w_n\leq v\}} \left(\frac{1}{b}\eta_n w_n-G(x,w_n)\right)
\end{equation}
and using $(g_1)$, $(g_2)$ we obtain the following uniform bound.
\begin{equation}\label{star 2}
\left|\int_{\{w_n\leq v\}} \left(\frac{1}{b}\eta_n w_n-G(x,w_n)\right)\right|\leq K\left[\left(\frac{1}{b}+1\right)v+\left(\frac{1}{b}+\frac{1}{r}\right)v^{r}\right]|\Omega|=C_1,~\forall n\in \mathbb{N}.
\end{equation}
 From Lemma $\ref{properties I}$ (3), we rewrite the last term of $\eqref{star}$ as follows.
	\begin{equation}\label{star 3}
	\int_\Omega \frac{1}{b}\theta_n w_n-F_\alpha(w_n)=\left(\frac{1}{b}-\frac{1}{p_s^*}\right)\int_{\{w_n>\alpha\}}|w_n|^{p_s^*}+\int_{\{w_n=\alpha\}}\frac{1}{b}\theta_n w_n+\frac{1}{p_s^*}\int_{\{w_n>\alpha\}}\alpha^{p_s^*}.
	\end{equation}
	We have used the fact that $F_\alpha(t)=\chi_{\{t\geq \alpha\}}\frac{1}{p_s^*}(|t|^{p_s^*}-\alpha^{p_s^*})$ for all $t\in \mathbb{R}$ to obatin $\eqref{star 3}$.\\
	On combining $\eqref{star}-\eqref{star 3}$ we get
	\begin{align}\label{star 4}
	J_\alpha(w_n)-\frac{1}{b}\langle v_n,w_n\rangle&\geq \left(\frac{1}{p}-\frac{1}{b}\right)\|w_n\|^p_{W_0^{s,p}(\Omega)}-\frac{\lambda}{1-\gamma} \|w_n\|_{L^{1-\gamma}(\Omega)}^{1-\gamma}- \mu C_1+\left(\frac{1}{b}-\frac{1}{p_s^*}\right)\int_{\{w_n>\alpha\}}|w_n|^{p_s^*}.
	\end{align}
	Using Theorem $\ref{embedding}$ and $\eqref{best}$ in $\eqref{star 4}$ we establish the following.
	\begin{align}\label{star 4-4}
		J_\alpha(w_n)-\frac{1}{b}\langle v_n,w_n\rangle&\geq \left(\frac{1}{p}-\frac{1}{b}\right)\|w_n\|^p_{W_0^{s,p}(\Omega)}-\frac{\lambda}{1-\gamma}|\Omega|^{\frac{p_s^*-1+\gamma}{p_s^*}}S_{s,p}^{-\frac{1-\gamma}{p}} \|w_n\|_{W_0^{s,p}(\Omega)}^{1-\gamma}- \mu C_1\nonumber\\&~~~~+\left(\frac{1}{b}-\frac{1}{p_s^*}\right)\int_{\{w_n>\alpha\}}|w_n|^{p_s^*}.
	\end{align}
	We already have $J_\alpha(w_n)=c+o_n(1)$ and $\|v_n\|_{(W_0^{s,p}(\Omega))^\prime}=o_n(1)$. Thus, 
	\begin{align}\label{star 5}
		J_\alpha(w_n)-\frac{1}{b}\langle v_n,w_n\rangle&\leq |J_\alpha(w_n)|+\frac{1}{b}\|v_n\|_{(W_0^{s,p}(\Omega))^\prime}\|w_n\|_{W_0^{s,p}(\Omega)}\nonumber\\&\leq c+1+\|w_n\|_{W_0^{s,p}(\Omega)}+o_n(1).
	\end{align}
	Since $p<b<p_s^*$ and $\alpha>0$, from $\eqref{star 4-4}$ and $\eqref{star 5}$, we obtain
	\begin{equation}\label{star 6}
	c+1+\|w_n\|_{W_0^{s,p}(\Omega)}+o_n(1)\geq \left(\frac{1}{p}-\frac{1}{b}\right)\|w_n\|^p_{W_0^{s,p}(\Omega)}-\frac{\lambda}{1-\gamma}|\Omega|^{\frac{p_s^*-1+\gamma}{p_s^*}}S_{s,p}^{-\frac{1-\gamma}{p}}\|w_n\|_{W_0^{s,p}(\Omega)}^{1-\gamma}- \mu C_1.
	\end{equation}
With the consideration of the above inequality $\eqref{star 6}$, we conclude that $(w_n)$ is a bounded sequence in $W_0^{s,p}(\Omega)$.
	\end{proof}
\end{lemma}
\begin{proposition}\label{ps}
	Let $(w_n)\subset W_0^{s,p}(\Omega)$ be a non-smooth (PS)$_c$ sequence such that $J_\alpha(w_n)\rightarrow c$ and $\Lambda_{J_\alpha}(w_n)\rightarrow 0$ with
	\begin{equation}\label{label}
	c<\left(\frac{1}{p}-\frac{1}{b}\right) S_{s,p}^{\frac{N}{sp}}-\left(\frac{1}{p}-\frac{1}{b}\right)^{-\frac{1-\gamma}{p-1+\gamma}} \left(\frac{\lambda}{1-\gamma}|\Omega|^{\frac{p_s^*-1+\gamma}{p_s^*}}S_{s,p}^{-\frac{1-\gamma}{p}}\right)^{\frac{p}{p-1+\gamma}}-\mu C_1=c^*,
	\end{equation} 
	where $S_{s,p}$ is defined in $\eqref{best}$ and $C_1$ is given in $\eqref{star 2}$. Then $J_\alpha$ satisfies the non-smooth (PS)$_c$ condition, i.e. $(w_n)$ admits a strongly convergent subsequence. Further, $c^*>0$ for a sufficiently small $\mu$.
	\begin{proof}
		According to Lemma $\ref{lemma 1}$, the (PS)$_c$ sequence $(w_n)$ is bounded in $W_0^{s,p}(\Omega)$. Let $v_n,\eta_n$ and $\theta_n$ are same as used in  the proof of Lemma $\ref{lemma 1}$. Then by $(g_2)$ and Lemma $\ref{properties I}$, we establish that $(\eta_n)$ and $(\theta_n)$ are bounded in $L^{\frac{r}{r-1}}(\Omega)$ and $L^{\frac{p_s^*}{p_s^*-1}}(\Omega)$, respectively. Hence, up to a subsequence,
		\begin{equation}\label{star 7}
		w_n\rightharpoonup w ~\text{in}~W_0^{s,p}(\Omega),~w_n(x)\rightarrow w(x)~\text{a.e. in}~ \Omega, \eta_n\overset{*}{\rightharpoonup}\eta ~\text{in}~L^{\frac{r}{r-1}}(\Omega)~\text{and}~\theta_n\overset{*}{\rightharpoonup}\theta ~\text{in}~L^{\frac{p_s^*}{p_s^*-1}}(\Omega).
		\end{equation}
		\begin{equation}\label{starr 1}
		\|w_n-w\|_{W_0^{s,p}(\Omega)}^p\rightarrow M.
		\end{equation}
		If $M=0$, then it implies $w_n\rightarrow w$ in $W_0^{s,p}(\Omega)$ as $n\rightarrow \infty$ and hence the proof. Thus, we assume $M>0$. Noting that $\eqref{star 7}$ indicates 
		$$\int_\Omega|w_n|^{p_s^*-1}\chi_{\{w_n>\alpha-\frac{1}{n}\}}\varphi dx\rightarrow \int_\Omega|w|^{p_s^*-1}\chi_{\{w>\alpha\}}\varphi dx,~\forall ~\varphi\in L^{p_s^*}(\Omega).$$
		Then, by Br\'{e}zis-Lieb Lemma (see \cite{Brezis}, Theorem 1), we have
		\begin{equation}\label{starr 2}
		\|w_n\|_{W_0^{s,p}(\Omega)}^p=\|w_n-w\|_{W_0^{s,p}(\Omega)}^p+\|w\|_{W_0^{s,p}(\Omega)}^p +o_n(1),
		\end{equation}
		\begin{equation}\label{starr 3}
		\int_\Omega|w_n\chi_{\{w_n>\alpha-\frac{1}{n}\}}|^{p_s^*}=\int_\Omega|w_n\chi_{\{w_n>\alpha-\frac{1}{n}\}}-w\chi_{\{w>\alpha\}}|^{p_s^*}+\int_\Omega|w\chi_{\{w_n>\alpha\}}|^{p_s^*}.
		\end{equation}
		Since $w_n\rightarrow w$ strongly in $L^q(\Omega)$ for any $1\leq q<p_s^*$ and $\int_{\Omega}\psi(x,w_n)<\int_\Omega\frac{1}{\underline{u}_\lambda^\gamma}$ (see the definition of $\psi$ given in $\eqref{defn}$), we obtain 
		\begin{equation}\label{starr 4}
		\left|\int_\Omega\psi(x,w_n)(w_n-w)\right|\leq o_n(1).
		\end{equation}
		Using Lemma $\ref{properties I}$, $(g_2)$ and $\eqref{starr 2}-\eqref{starr 4}$ we obtain
		\begin{align}
		o_n(1)&=\langle v_n,w_n-w\rangle\nonumber\\&=\|w_n\|^p_{W_0^{s,p}(\Omega)}-\int_{\mathbb{R}^{2N}}\frac{|w_n(x)-w_n(y)|^{p-2}(w_n(x)-w_n(y))}{|x-y|^{N+2s}}(w(x)-w(y))dxdy\nonumber \\&~~~~-\int_\Omega\psi(x,w_n)(w_n-w)-\int_\Omega \theta_n (w_n-w) +o_n(1)\nonumber\\&=\|w_n-w\|_{W_0^{s,p}(\Omega)}^p-\int_\Omega |w_n|^{p_s^*-1}(w_n-w)\chi_{\{w_n>\alpha-\frac{1}{n}\}}dx+o_n(1)\nonumber\\&=\|w_n-w\|_{W_0^{s,p}(\Omega)}^p-\int_\Omega |w_n|^{p_s^*}\chi_{\{w_n>\alpha-\frac{1}{n}\}}dx+\int_\Omega |w|^{p_s^*}\chi_{\{w>\alpha\}}dx+o_n(1)\nonumber\\&= \|w_n-w\|_{W_0^{s,p}(\Omega)}^p-\int_\Omega \left|w_n\chi_{\{w_n>\alpha-\frac{1}{n}\}}-w\chi_{\{w>\alpha\}}\right|^{p_s^*}dx+o_n(1).\nonumber
		\end{align}
		This implies
		\begin{equation}\label{starr 5}
		\|w_n-w\|_{W_0^{s,p}(\Omega)}^p= \int_\Omega \left|w_n\chi_{\{w_n>\alpha-\frac{1}{n}\}}-w\chi_{\{w>\alpha\}}\right|^{p_s^*}dx+o_n(1).
		\end{equation}
		From $\eqref{best}$, $\eqref{starr 1}$, $\eqref{starr 2}$ and $\eqref{starr 5}$ we have
		\begin{align}\label{starr 6}
		S_{s,p}&\leq \frac{	\|w_n-w\|_{W_0^{s,p}(\Omega)}^p}{\left(\int_\Omega \left|w_n\chi_{\{w_n>\alpha-\frac{1}{n}\}}-w\chi_{\{w>\alpha\}}\right|^{p_s^*}dx\right)^{p/p^*_s}}\nonumber\\&\leq M^{\frac{sp}{N}}+o_n(1)
		\end{align}
		and hence $M>S_{s,p}^{\frac{N}{sp}}+o_n(1)$. According to $\eqref{star 4}$ and the fact that $p<b<p_s^*$, we get
			\begin{align}\label{starr 7}
			J_\alpha(w_n)-\frac{1}{b}\langle v_n,w_n\rangle&\geq \left(\frac{1}{p}-\frac{1}{b}\right)\|w_n\|^p_{W_0^{s,p}(\Omega)}-\frac{\lambda}{1-\gamma} \|w_n\|_{L^{1-\gamma}(\Omega)}^{1-\gamma}- \mu C_1,
			\end{align}
			where $b$ is given in $(g_3)$ and $C_1$ is obtained in $\eqref{star 2}$ which is independent of $\alpha$. Since $J_\alpha(w_n)=c+o_n(1)$ and $\|v_n\|_{(W_0^{s,p}(\Omega))^\prime}=o_n(1)$, using $\eqref{star 7}-\eqref{starr 2}$, $\eqref{starr 6}$, $\eqref{starr 7}$, Theorem $\ref{embedding}$, $\eqref{best}$ and Young's inequality we obtain
			\begin{align}\label{starr 8}
			c&\geq \left(\frac{1}{p}-\frac{1}{b}\right)(M+\|w\|_{W_0^{s,p}(\Omega)}^p)-\frac{\lambda}{1-\gamma}|\Omega|^{\frac{p_s^*-1+\gamma}{p_s^*}}S_{s,p}^{-\frac{1-\gamma}{p}} \|w\|_{W_0^{s,p}(\Omega)}^{1-\gamma}- \mu C_1+o_n(1)\nonumber\\&\geq \left(\frac{1}{p}-\frac{1}{b}\right)(M+\|w\|_{W_0^{s,p}(\Omega)}^p)- \left(\frac{1}{p}-\frac{1}{b}\right) \|w\|_{W_0^{s,p}(\Omega)}^p\nonumber\\&~~~~-\left(\frac{1}{p}-\frac{1}{b}\right)^{-\frac{1-\gamma}{p-1+\gamma}} \left(\frac{\lambda}{1-\gamma}|\Omega|^{\frac{p_s^*-1+\gamma}{p_s^*}}S_{s,p}^{-\frac{1-\gamma}{p}}\right)^{\frac{p}{p-1+\gamma}}-\mu C_1+o_n(1)\nonumber\\&\geq \left(\frac{1}{p}-\frac{1}{b}\right) S_{s,p}^{\frac{N}{sp}}-\left(\frac{1}{p}-\frac{1}{b}\right)^{-\frac{1-\gamma}{p-1+\gamma}} \left(\frac{\lambda}{1-\gamma}|\Omega|^{\frac{p_s^*-1+\gamma}{p_s^*}}S_{s,p}^{-\frac{1-\gamma}{p}}\right)^{\frac{p}{p-1+\gamma}}-\mu C_1+o_n(1)\nonumber\\
			&=c^*+o_n(1).
			\end{align}
			The above inequality $\eqref{starr 8}$ is a contradiction to $\eqref{label}$. Thus, $M=0$ and $w_n\rightarrow w$ strongly in $W_0^{s,p}(\Omega)$. Moreover, $c^*>0$ for a sufficiently small $\mu>0$.
				\end{proof}
			\end{proposition}
	\noindent In the next lemma, we prove the functional $J_\alpha$ satisfies all the hypotheses of Theorem $\ref{mountain pass non smooth}$ with a suitable choice of $\beta$ given in $(g_4)$.
	\begin{lemma}\label{lemma 2}
		Let $(g_1)-(g_5)$ hold. Then there exist $\rho_1>0$, $\bar{\mu}\in(0,1)$, $\bar{\lambda}=\bar{\lambda}(\rho_1)>0$, $\rho_2>0$, $m^*>0$ and $\sigma\in W_0^{s,p}(\Omega)$ such that for every $\alpha>0$, every $\lambda\in (0,\bar{\lambda}]$ and every $\mu\in(0,\bar{\mu}]$ we have
		\begin{enumerate}
			\item $\underset{m\in [0,m^*]}{\sup}J_\alpha(m\sigma)<\left(\frac{1}{p}-\frac{1}{b}\right) S_{s,p}^{\frac{N}{sp}}-\left(\frac{1}{p}-\frac{1}{b}\right)^{-\frac{1-\gamma}{p-1+\gamma}} \left(\frac{\lambda}{1-\gamma}|\Omega|^{\frac{p_s^*-1+\gamma}{p_s^*}}S_{s,p}^{-\frac{1-\gamma}{p}}\right)^{\frac{p}{p-1+\gamma}}-\mu C_1$ where $S_{s,p}$ is defined in $\eqref{best}$ and $C_1$ is given in $\eqref{star 2}$.
			\item $J_\alpha(w)\geq \rho_2$ for every $w\in W_0^{s,p}(\Omega)$, $\|w\|_{W_0^{s,p}(\Omega)}=\rho_1$, where $\rho_1$ and $\rho_2$ are independent of $\alpha$.
			\item $J_\alpha(m^*\sigma)<0$ for $\|m^*\sigma\|_{W_0^{s,p}(\Omega)}> \rho_{1}$.
		\end{enumerate}
		\begin{proof}
		Let us fix $\sigma\in W_0^{s,p}(\Omega)$ with $\sigma>0$ in $\Omega$ and $\|\sigma\|_{W_0^{s,p}(\Omega)}=1$. By the hypothesis $(g_4)$ we can write
		\begin{equation}\label{starrr 1}
		J_\alpha(m\sigma)=\frac{1}{p}m^p-m\int_{\{m\sigma>\beta\}}\sigma dx +\beta\left|\{m\sigma>\beta\}\right|,~\forall ~m\geq 0.
		\end{equation}
	It is easy to check that there exist $\bar{\mu}\in (0,1)$ and $\lambda^*>0$ such that for any $\lambda\in (0,\lambda^*]$ and $\mu\in(0,\bar{\mu}]$ we have $$\left(\frac{1}{p}-\frac{1}{b}\right) S_{s,p}^{\frac{N}{sp}}-\left(\frac{1}{p}-\frac{1}{b}\right)^{-\frac{1-\gamma}{p-1+\gamma}} \left(\frac{\lambda}{1-\gamma}|\Omega|^{\frac{p_s^*-1+\gamma}{p_s^*}}S_{s,p}^{-\frac{1-\gamma}{p}}\right)^{\frac{p}{p-1+\gamma}}-{\mu} C_1>0.$$
Hence, we choose $m^*>0$ such that
\begin{align}\label{starrr 2}
\frac{1}{p}(m^*)^p-m^*\int_{\{m\sigma>\beta\}}\sigma dx<0
\end{align} and
\begin{align}\label{starrr 3}
\frac{(m^*)^p}{p}<\left(\frac{1}{p}-\frac{1}{b}\right) S_{s,p}^{\frac{N}{sp}}-\left(\frac{1}{p}-\frac{1}{b}\right)^{-\frac{1-\gamma}{p-1+\gamma}} \left(\frac{\lambda}{1-\gamma}|\Omega|^{\frac{p_s^*-1+\gamma}{p_s^*}}S_{s,p}^{-\frac{1-\gamma}{p}}\right)^{\frac{p}{p-1+\gamma}}-{\mu} C_1.
\end{align}	
Since $m^*$ does not depend on $\beta$, we get
\begin{equation}\label{starrr 4}
\int_\Omega\sigma dx=\int_{\{m^*\sigma>\beta\}}\sigma dx+o_\beta(1).
\end{equation}
Combining $\eqref{starrr 1}-\eqref{starrr 4}$ we choose $\beta>0$ very small such that $J_\alpha(m^*\sigma)<0$ for every $\sigma>0$. This proves $(3)$ and also $(1)$.\\
According to $(g_1)$, $(g_5)$ and the fact that $J_\alpha(0)=0$, there exist $\overline{K}$ and $c_{\overline{K}}>0$ (independent of $\alpha$) such that $\overline{K}<\lambda_1$ and 
$$\left|g(x,t)\right|\leq \overline{K}|t|^{p-1}+c_{\overline{K}}|t|^{r-1}, ~\forall~ t\in\mathbb{R}.$$
Thus, from $\eqref{first}$ for every $w\in W_0^{s,p}(\Omega)$
\begin{equation}\label{starrr 5}
\int_\Omega G(x,w) dx\leq \frac{\overline{K}}{p\lambda_1}\iint_{\mathbb{R}^{2N}}\frac{|w(x)-w(y)|^p}{|x-y|^{N+sp}}dxdy+\frac{c_{\overline{K}}}{r}\int_\Omega |w|^r dx.
\end{equation}
Hence,
\begin{equation}\label{starrr 6}
J_\alpha(w)\geq \frac{1}{p}\left(1-\frac{\mu\overline{K}}{\lambda_1}\right)\|w\|^p_{W_0^{s,p}(\Omega)}-C_1 c_{\overline{K}}\mu \|w\|^r_{W_0^{s,p}(\Omega)}-C_2\|w\|^{p_s^*}_{W_0^{s,p}(\Omega)}-\frac{\lambda}{1-\gamma}|\Omega|^{\frac{p_s^*-1+\gamma}{p_s^*}}S_{s,p}^{-\frac{1-\gamma}{p}}\|w\|^{1-\gamma}_{W_0^{s,p}(\Omega)},
\end{equation}
for some constant $C_1,C_2>0$ independent of $\alpha$. Since $1-\gamma<1<p<r<p_s^*$, the function 
$$h(t)=\frac{1}{p}\left(1-\frac{\mu\overline{K}}{\lambda_1}\right)t^{p-1+\gamma}-C_1 c_{\overline{K}}\mu t^{r-1+\gamma}-C_2t^{p_s^*-1+\gamma}, ~t\in[0,1]$$
admits a maximum at some $\rho_1\in (0,1]$ small enough, i.e. $\underset{t\in[0,1]}{\max}h(t)=h(\rho_1)>0$. Therefore, let 
$$\lambda^{**}=\frac{(1-\gamma)S_{s,p}^{\frac{1-\gamma}{p}}}{2|\Omega|^{\frac{p_s^*-1+\gamma}{p_s^*}}}h(\rho_1),$$
then for every $w\in W_0^{s,p}(\Omega)$ with $\|w\|_{W_0^{s,p}(\Omega)}=\rho_1\leq 1$ and for every $\lambda\in (0,\lambda^{**}]$, we have $J_\alpha(w)\geq \rho_1^{1-\gamma}h(\rho_1)/2=\rho_2$. Let $\bar{\lambda}=\min\{\lambda^*,\lambda^{**}\}$. Then with $\lambda\in (0,\bar{\lambda})$ and $\mu\in (0,\bar{\mu})$, we conclude the proof.
		\end{proof}
	\end{lemma}
	\begin{proposition}\label{proof part 1}
		Let $(g_1)-(g_5)$ are satisfied. Then there exist $\bar{\alpha},\bar{\lambda}, \bar{\mu}
		>0$ such that for every $a\in (0,\bar{\alpha})$, every $\lambda\in (0,\bar{\lambda})$ and every $\mu\in (0,\bar{\mu})$ the problem $\eqref{main prob intro}$ admits at least one nontrivial weak solution $u_\alpha$. Furthermore, the lebesgue measure of the set $\{x\in\Omega:u_\alpha>\alpha\}$ is positive.
		\begin{proof}
			Let  $$c_\alpha=\underset{\zeta\in \Gamma}{\inf}~\underset{t\in [0,1]}{\max}~J_\alpha(\zeta(t))~\text{and}~\Gamma=\{\zeta\in C([0,1];W_0^{s,p}(\Omega)):\zeta(0)=0,\zeta(1)=m^*\sigma\},$$
		where $m^*$, $\sigma$, $\rho_1$, $\rho_2$, $\bar{\lambda}$, $\bar{\mu}$ are obtained in Lemma $\ref{lemma 2}$. Since $J_\alpha$ satisfy the hypothese of Theorem $\ref{mountain pass non smooth}$ (refer Lemma $\ref{lemma 2}$), we guarantee the existence of a sequence $(w_n)\subset W_0^{s,p}(\Omega)$ that satisfy $J_\alpha(w_n)=c_\alpha+o_n(1)$ and $\Lambda_{J_\alpha}(w_n)=o_n(1)$. By $(1)$ and $(2)$ of Lemma $\ref{lemma 2}$, we also have 
		\begin{equation}\label{starrr 7}
		\rho_2\leq c_\alpha<\left(\frac{1}{p}-\frac{1}{b}\right) S_{s,p}^{\frac{N}{sp}}-\left(\frac{1}{p}-\frac{1}{b}\right)^{-\frac{1-\gamma}{p-1+\gamma}} \left(\frac{\lambda}{1-\gamma}|\Omega|^{\frac{p_s^*-1+\gamma}{p_s^*}}S_{s,p}^{-\frac{1-\gamma}{p}}\right)^{\frac{p}{p-1+\gamma}}-{\mu} C_1, ~\forall ~\alpha>0.
		\end{equation}
		From Proposition $\ref{ps}$ there exists $w_\alpha\in W_0^{s,p}(\Omega)$ such that, up to a subsequence, $w_n\rightarrow w_\alpha$ in $W_0^{s,p}(\Omega)$ as $n\rightarrow \infty$. This implies $J_\alpha(w_\alpha)=c_\alpha$ and $0\in \partial J_\alpha(w_\alpha)$. Thus, by $\eqref{differential}$ and Lemma $\ref{properties I}$, there exist $\eta_\alpha\in L^{\frac{r}{r-1}}(\Omega)$ and $\theta_\alpha\in L^{\frac{p_s^*}{p_s^*-1}}(\Omega)$ such that 
		\begin{equation}\label{starrr 8}
		\int_{\mathbb{R}^{2N}}\frac{|w_\alpha(x)-w_\alpha(y)|^{p-2}(w_\alpha(x)-w_\alpha(y)}{|x-y|^{N+sp}}(\varphi(x)-\varphi(y))dxdy=\lambda\int_\Omega \psi(x,w_\alpha)\varphi dx+\mu\langle\eta_\alpha,\varphi\rangle+\langle\theta_\alpha,\varphi\rangle,
		\end{equation}
		for every $\varphi\in W_0^{s,p}(\Omega)$, where $\eta_\alpha(x)\in [\underline{g}(x,w_\alpha(x)),\overline{g}(x,w_\alpha(x))]\geq 0$ and $\theta_\alpha(x)\in [\underline{f}(w_\alpha(x)),\overline{g}(w_\alpha(x))]\geq 0$ a.e. in $\Omega$. According to the strong maximum principle (Lemma 2.3, \cite{Mosconi}) we have $w_\alpha>0$ a.e. in $\Omega$. \\
		This proves that $w_\alpha$ is a weak solution to $\eqref{main prob help}$. By the weak comparison principle for fractional $p$-Laplacian (Lemma 3.1, \cite{Ghosh}), we conclude that $\underline{u}_\lambda\leq w_\alpha$ a.e. in $\Omega$. This implies $\psi(x,w_\alpha)=w_\alpha^{-\gamma}$ a.e. in $\Omega$ and $w_\alpha=u_\alpha$ is a weak solution to $\eqref{main prob intro}$.\\
		The next claim is to prove that the set $\{x\in \Omega:w_\alpha(x)>\alpha\}$ has positive lebesgue measure in $\mathbb{R}^N$. We prove this claim by method of contradiction. For this let us assume that the set $\{x\in \Omega:w_\alpha(x)>\alpha\}$ is of zero lebesgue measure in $\mathbb{R}^N$. Thus, $w_\alpha(x)\leq \alpha$ a.e. in $\Omega$.\\
		From Lemma $\ref{properties I}$, $\eqref{starrr 5}$ and $\eqref{starrr 8}$ we obtain 
		\begin{align}
	\|w_\alpha\|_{W_0^{s,p}(\Omega)}^p&=\mu\int_{\Omega}\eta_\alpha w_\alpha+\int_{\{w_\alpha=\alpha\}} \theta_\alpha w_\alpha+\lambda\int_\Omega \psi(x,w_\alpha)w_\alpha\nonumber\\&\leq \mu\int_\Omega(\overline{K}\alpha^p+c_{\overline{K}}\alpha^r)dx+\int_{\{w_\alpha=\alpha\}}\alpha^{p_s^*}+\lambda\int_\Omega\alpha^{1-\gamma}.\nonumber
		\end{align}
	Since $J_\alpha(w_\alpha)=c_\alpha$ by Lemma $\ref{lemma 2}$ and $\eqref{starrr 7}$, for $\alpha>0$ small enough, we establish
	$$p\rho_2\leq [\mu(\overline{K}+c_{\overline{K}})+1+\lambda]|\Omega|\alpha^{1-\gamma}.$$
This contradicts the fact that $\rho_2$ is independent of $\alpha$. Thus, there exists $\bar{\alpha}>0$ small such that for any $\alpha\in (0,\bar{\alpha})$ the set $\{x\in \Omega:w_\alpha(x)>\alpha\}$ has positive lebesgue measure. 
		\end{proof}
	\end{proposition}
		\section{Proof of Theorem $\ref{main result heaviside}$ (2)}\label{4}
\noindent Let $u_\alpha$ be a nontrivial weak solution to $\eqref{main prob intro}$ given in Proposition $\ref{proof part 1}$. In this section, we prove the second part of Theorem $\ref{main result heaviside}$, i.e. we examine the nature of $(u_\alpha)$ as $\alpha\rightarrow 0^+$.\\
Consider the functional $J_0^0:W_0^{s,p}(\Omega)\rightarrow\mathbb{R}$ associated to $\eqref{main prob 0}$ defined by
$$J_0^0(u)=\frac{1}{p}\int_{\mathbb{R}^{2N}}\frac{|u(x)-u(y)|^{p}}{|x-y|^{N+sp}}dxdy-\frac{\lambda}{1-\gamma}\int_\Omega u^{1-\gamma}-\mu\int_\Omega G(x,u)-\frac{1}{p_s^*}\int_\Omega (u^+)^{p_s^*},~\forall ~u\in W_0^{s,p}(\Omega).$$
Let us define
\begin{equation}\label{starrrr}
c_0=\underset{\zeta\in \Gamma}{\inf}~\underset{t\in [0,1]}{\max}~J_0^0(\zeta(t))~\text{and}~\Gamma=\{\zeta\in C([0,1];W_0^{s,p}(\Omega)):\zeta(0)=0,\zeta(1)=m^*\sigma\},
\end{equation}
where $m^*,\sigma$ as obtained in Lemma $\ref{lemma 2}$.
\begin{lemma}\label{lemma 3}
	$\lim\limits_{\alpha\rightarrow 0^+}c_\alpha=c_0\geq \rho_2$, where $c_\alpha$, $c_0$ and $\rho_2$ are given in Proposition $\ref{proof part 1}$, $\eqref{starrrr}$ and Lemma $\ref{lemma 2}$, respectively.
	\begin{proof}
		Since $F_\alpha(t)=\chi_{\{t\geq \alpha\}}\frac{1}{p_s^*}(|t|^{p_s^*}-\alpha^{p_s^*})$, we obtain
		\begin{align}\label{starrrr 1}
		\left|\frac{1}{p_s^*}\int_\Omega(u^+)^{p_s^*}-\int_\Omega F_\alpha(u)\right|&=	\left|\frac{1}{p_s^*}\int_\Omega(u^+)^{p_s^*}\chi_{\{u\leq \alpha\}}+\frac{1}{p_s^*}\int_\Omega \alpha^{p_s^*}\chi_{\{u>\alpha\}}\right|\nonumber\\&\leq \frac{2\alpha^{p_s^*}|\Omega|}{p_s^*}.
		\end{align}
Clearly, $J_0^0(u)\leq J_\alpha^0(u)$, for all $u\in W_0^{s,p}(\Omega)$. Thus, $c_0\leq c_\alpha$, for all $\alpha>0$. According to $\eqref{starrrr 1}$ we establish
\begin{equation}\label{starrrr 2}
J_\alpha^0(u)=J_0^0(u)+o_\alpha(1),~\forall~u\in W_0^{s,p}(\Omega), 
\end{equation}
where $o_\alpha(1)\rightarrow0$ as $\alpha\rightarrow0^+$ independently of $u$. This gives, 
\begin{equation}\label{starrrr 3}
J_\alpha^0(\zeta(t))=J_0^0(\zeta(t))+o_\alpha(1),~\forall~\zeta\in \Gamma,~t\in[0,1]
\end{equation}
and hence $c_\alpha=c_0+o_\alpha(1)$.\\
With the consideration of $\eqref{starrrr 3}$ and Lemma $\ref{lemma 2}$ $(3)$, we conclude that 	$\lim\limits_{\alpha\rightarrow 0^+}c_\alpha=c_0\geq \rho_2$. 
	\end{proof}
	\begin{proposition}
		For any sequence $\alpha_n\in (0,\bar{\alpha})$ with $\alpha_n\rightarrow 0^+$, we have, up to a subsequence, $u_{\alpha_n}\rightarrow u_0$ in $W_0^{s,p}(\Omega)$, where $u_0$ is a nontrivial weak solution to the problem $\eqref{main prob 0}$.
		\begin{proof}
			Let $u_\alpha$ be the weak solutionto $\eqref{main prob intro}$ given in Proposition $\ref{proof part 1}$. Thus, $J_\alpha^0(u_\alpha)=J_\alpha(u_\alpha)=c_\alpha$and  
			\begin{equation}\label{starrrr 4}
			\int_{\mathbb{R}^{2N}}\frac{|u_\alpha(x)-u_\alpha(y)|^{p-2}(u_\alpha(x)-u_\alpha(y)}{|x-y|^{N+sp}}(\varphi(x)-\varphi(y))dxdy-\mu\int_\Omega \eta_\alpha \varphi dx-\int_\Omega\theta_\alpha\varphi dx-\lambda\int_\Omega\frac{\varphi}{u_\alpha^{\gamma}}dx=0,
			\end{equation}
			where $\eta_\alpha(x)\in [\underline{g}(x,u_\alpha(x)),\overline{g}(x,u_\alpha(x))]$ and $\theta_\alpha(x)\in [\underline{f}_\alpha(u_\alpha(x)),\overline{f}_\alpha(u_\alpha(x))]$ a.e. in $\Omega$ with $f_\alpha(t)=H(t-\alpha)t^{p_s^*-1}$.\\
	Consider the sequence $(w_n)\subset W_0^{s,p}(\Omega)$ obtained in Proposition $\ref{proof part 1}$ with $w_n\rightarrow u_\alpha$ in $W_0^{s,p}(\Omega)$. From $\eqref{star 6}$ we have
		\begin{equation}
		c_\alpha+1+\|w_n\|_{W_0^{s,p}(\Omega)}\geq \left(\frac{1}{p}-\frac{1}{b}\right)\|w_n\|^p_{W_0^{s,p}(\Omega)}-\frac{\lambda}{1-\gamma}|\Omega|^{\frac{p_s^*-1+\gamma}{p_s^*}}S_{s,p}^{-\frac{1-\gamma}{p}}\|w_n\|_{W_0^{s,p}(\Omega)}^{1-\gamma}- \mu C_1,~\forall ~\alpha>0,\nonumber
		\end{equation}
		where $C_1$ is independent of $\alpha$ (refer $\eqref{star 2}$). Thus,
		\begin{equation}
		c_\alpha+1+\|w_\alpha\|_{W_0^{s,p}(\Omega)}\geq \left(\frac{1}{p}-\frac{1}{b}\right)\|w_\alpha\|^p_{W_0^{s,p}(\Omega)}-\frac{\lambda}{1-\gamma}|\Omega|^{\frac{p_s^*-1+\gamma}{p_s^*}}S_{s,p}^{-\frac{1-\gamma}{p}}\|w_\alpha\|_{W_0^{s,p}(\Omega)}^{1-\gamma}- \mu C_1,~\forall ~\alpha>0\nonumber
		\end{equation}
		and the sequence $(u_\alpha)$ is uniformly bounded in $W_0^{s,p}(\Omega)$. By $(g_2)$ and Lemma $\ref{properties I}$, we establish that $(\eta_\alpha)$ and $(\theta_\alpha)$ are bounded in $L^{\frac{r}{r-1}}(\Omega)$ and $L^{\frac{p_s^*}{p_s^*-1}}(\Omega)$, respectively.\\ Consider a sequence $(\alpha_n)\subset(0,\bar{\alpha})$ with $\alpha_n\rightarrow 0^+$. Hence, up to a subsequence,
		\begin{align}\label{starrrr 5}
		u_{\alpha_n}\rightharpoonup u_0 ~\text{in}~W_0^{s,p}(\Omega),&~	u_{\alpha_n}(x)\rightarrow 	u_0(x)~\text{a.e. in}~ \Omega,~ u_{\alpha_n}\rightarrow u_\alpha~ \text{in}~ L^q(\Omega) ~\text{for any}~ 1\leq q<p_s^*,	\nonumber\\ \eta_{\alpha_n}\overset{*}{\rightharpoonup}\eta_0 ~\text{in}~L^{\frac{r}{r-1}}(\Omega)~&\text{and}~\theta_{\alpha_n}\overset{*}{\rightharpoonup} u_0^{p_s^*-1} ~\text{in}~L^{\frac{p_s^*}{p_s^*-1}}(\Omega).
		\end{align}
We already have $u_{\alpha_n}\geq \underline{u}_\lambda$ for all $n\in \mathbb{N}$ (refer Proposition $\ref{proof part 1}$), where $\underline{u}_\lambda$ is a weak solution to $\eqref{main prob help}$ given in Lemma $\ref{main result help}$. Thus, by combining $\eqref{starrrr 4}$ and $\eqref{starrrr 5}$ we pass the limit $\alpha_n\rightarrow 0^+$ to obtain
	\begin{equation}
	\int_{\mathbb{R}^{2N}}\frac{|u_0(x)-u_0(y)|^{p-2}(u_0(x)-u_0(y)}{|x-y|^{N+sp}}(\varphi(x)-\varphi(y))dxdy-\mu\int_\Omega \eta_0 \varphi dx-\int_\Omega u_0^{p_s^*-1}\varphi dx-\lambda\int_\Omega\frac{\varphi}{u_0^{\gamma}}dx=0,\nonumber
	\end{equation}
for every $\varphi\in W_0^{s,p}(\Omega)$. From Proposition $\ref{convergence}$, $\eta_0\in \partial I_\mu(u_0)$. \\
According to $\eqref{starrr 7}$ we have
\begin{equation}\label{starrrr 6}
\rho_2\leq J^0_{\alpha_n}(u_{\alpha_n})=c_\alpha<\left(\frac{1}{p}-\frac{1}{b}\right) S_{s,p}^{\frac{N}{sp}}-\left(\frac{1}{p}-\frac{1}{b}\right)^{-\frac{1-\gamma}{p-1+\gamma}} \left(\frac{\lambda}{1-\gamma}|\Omega|^{\frac{p_s^*-1+\gamma}{p_s^*}}S_{s,p}^{-\frac{1-\gamma}{p}}\right)^{\frac{p}{p-1+\gamma}}-{\mu} C_1, ~\forall ~\alpha_n>0
\end{equation}
and from $\eqref{starrrr 2}$, 
\begin{equation}\label{starrrr 7}
c_{\alpha_n}=J^0_{\alpha_n}(u_{\alpha_n})=J^0_{0}(u_{\alpha_n})+o_{\alpha_n}(1).
\end{equation}
Considering $\eqref{starrrr 6}$, $\eqref{starrrr 7}$ and then following the proof of Proposition $\ref{ps}$, we get
\begin{equation}\label{starrrr 8}
u_{\alpha_n}\rightarrow u_0~\text{ in }W_0^{s,p}(\Omega)~\text{as}~\alpha_n\rightarrow0^+.
\end{equation}
Combining $\eqref{starrrr 8}$, Lemma $\ref{lemma 3}$ and $\eqref{starrrr 2}$, we have $J_0^0(u_0)=c_0\geq \rho_2>0$. This proves that $u_0$ is a nontrivial weak solution to $\eqref{main prob 0}$ and we conclude the proof.
		\end{proof}
	\end{proposition}
\end{lemma}


\section*{Acknowledgement}
\noindent The author Debajyoti Choudhuri thanks the grant received from Council of Scientific and Industrial Research (CSIR), India for the research grant (09/983(0013)/2017-EMR-I). The author Akasmika Panda thanks the financial assistantship received from the Ministry of Human Resource Development (M.H.R.D.), Govt. of India. 

\section{References}


\end{document}